\magnification=1200

\catcode`\@=11
 
\def\hexnumber@#1{\ifnum#1<10 \number#1\else
 \ifnum#1=10 A\else\ifnum#1=11 B\else\ifnum#1=12 C\else
 \ifnum#1=13 D\else\ifnum#1=14 E\else\ifnum#1=15 F\fi\fi\fi\fi\fi\fi\fi}
 
\font\tenmsb=msbm10
\font\sevenmsb=msbm7
\font\fivemsb=msbm5
\newfam\msbfam
\textfont\msbfam=\tenmsb  \scriptfont\msbfam=\sevenmsb
  \scriptscriptfont\msbfam=\fivemsb
\def\msb@{\hexnumber@\msbfam}

\mathchardef\nodiv="3\msb@2D
\mathchardef\varemptyset="0\msb@3F

\font\tenbbb=bbold10
\font\sevenbbb=bbold7
\font\fivebbb=bbold5
\newfam\bbbfam
\textfont\bbbfam=\tenbbb  \scriptfont\bbbfam=\sevenbbb
  \scriptscriptfont\bbbfam=\fivebbb
\def\Bbb{\relax\ifmmode\let\next\Bbb@\else
 \def\next{\errmessage{Use \string\Bbb\space only in math mode}}\fi\next}
\def\Bbb@#1{{\Bbb@@{#1}}}

\def\Bbb@@#1{\fam\msbfam#1}    
 
\font\tenscr=rsfs10 
\font\sevenscr=rsfs7 
\font\fivescr=rsfs5 
\skewchar\tenscr='177 \skewchar\sevenscr='177 \skewchar\fivescr='177
\newfam\scrfam
\textfont\scrfam=\tenscr \scriptfont\scrfam=\sevenscr
  \scriptscriptfont\scrfam=\fivescr
\def\Scr{\relax\ifmmode\let\next\Scr@\else
 \def\next{\errmessage{Use \string\Scr\space only in math mode}}\fi\next}
\def\Scr@#1{{\Scr@@{#1}}}
\def\Scr@@#1{\fam\scrfam#1}

\catcode`\@=12

\def\R{{\Bbb R}}

\def\defeq{\buildrel \rm def \over =}

\newcount\thmno \thmno=1
\long\def\theorem#1#2{\edef#1{Theorem~\the\thmno}\noindent
    {\sc Theorem \the\thmno.\enspace}{\it #2}\endgraf
    \penalty55\global\advance\thmno by 1}
\newcount\lemmano \lemmano=1
\long\def\lemma#1#2{\edef#1{Lemma~\the\lemmano}\noindent
    {\sc Lemma \the\lemmano.\enspace}{\it #2}\endgraf
    \penalty55\global\advance\lemmano by 1}
\newcount\propno \propno=1
\long\def\proposition#1#2{\edef#1{Proposition~\the\propno}\noindent
    {\sc Proposition \the\propno.\enspace}{\it #2}\endgraf
    \penalty55\global\advance\propno by 1}
\newcount\corolno \corolno=1
\long\def\corollary#1#2{\edef#1{Corollary~\the\corolno}\noindent
    {\sc Corollary \the\corolno.\enspace}{\it #2}\endgraf
    \penalty55\global\advance\corolno by 1}

\newcount\conjno \conjno=1
\long\def\conjecture#1#2{\edef#1{Conjecture~\the\conjno}\noindent
    {\sc Conjecture \the\conjno.\enspace}{\it #2}\endgraf
    \global\advance\conjno by 1}

\newcount\eqnno \eqnno=0
\def\eqn#1{\global\edef#1{(\the\secno.\the\eqnno)}#1\global
   \advance\eqnno by 1}
\newcount\secno \secno=0
\def\sec#1\par{\eqnno=1\global\advance\secno by 1\bigskip
    \noindent\centerline{\sc\the\secno. #1}\par\nobreak\noindent}

\def\definition{\noindent{\sc Definition.\enspace}}
\def\remark{\noindent{\sc Remark.\enspace}}
\def\proof{\noindent{\it Proof.\enspace}}
\def\qed{~\vrule width4pt height6pt depth1pt}
\font\sc=cmcsc10

\parskip=12pt plus4pt minus4pt
\def\itc#1{{\it #1\/}}

\def\I{{\Scr I}}

\def\SB{{\Scr B}}

\font\smrm=cmr9
\font\smbf=cmbx9

\font\smit=cmti9
\font\smsc=cmcsc9
\font\bigbf=cmbx12 

\vbox to \vsize{
\vfill
\centerline{\bigbf An Independent Set Axiomatization of Symplectic Matroids}
\bigskip
\centerline{Timothy Y. Chow}
\centerline{University of Michigan}
\centerline{Department of Mathematics}
\centerline{Ann Arbor, MI 48109-1109}
\centerline{\tt tchow@umich.edu}
\vfill\vfill}

\sec Introduction

Symplectic matroids are rather esoteric objects,
so they require more than the average amount of motivation.
Thus we shall give a broader introduction than is strictly necessary
for understanding our main result.

Symplectic matroids occur at the confluence of two streams
in modern combinatorics: Coxeter groups and matroid theory.
Recall that a \itc{Coxeter system} is a pair $(W,S)$
where $W$ is a group and $S\subseteq W$ is a set of generators
subject only to relations of the form
$$(ss')^{m(s,s')} = 1,$$
where $m(s,s)=1$ and $m(s,s')=m(s',s)\ge 2$ if $s\ne s'$.
The group $W$ is called a \itc{Coxeter group}.

Perhaps the most famous examples of Coxeter groups are
those that arise in Lie theory.
There are exactly four infinite families
of finite-dimensional complex simple Lie algebras
(plus five exceptions that we will ignore here):
$$\eqalign{A_n&: \hbox{special linear Lie algebras}\cr
 B_n&: \hbox{odd orthogonal Lie algebras}\cr
 C_n&: \hbox{symplectic Lie algebras}\cr
 D_n&: \hbox{even orthogonal Lie algebras}\cr}$$
To each of these Lie algebras is associated
a finite group called its \itc{Weyl group}.
The precise definition of a Weyl group is not important here;
what is important is that it is always a Coxeter group.
For example, the Weyl group of~$A_n$ is
the symmetric group on $n+1$ letters,
and it is not too difficult to verify directly that
the symmetric group is a Coxeter group
(let $S$ be the set of all adjacent transpositions $(i, i+1)$).

Coxeter groups typically arise in combinatorics when
some combinatorial concept is shown
to be definable purely in terms of the symmetric group.
Since the symmetric group is a Coxeter group,
one can try replacing the term ``symmetric group'' in such a definition
by an arbitrary Coxeter group
to see if the result makes sense and is interesting.
This simple tactic has turned out to be surprisingly fruitful,
in part because it often reveals how
techniques from Lie theory and other areas of mathematics
can be brought to bear on combinatorial problems.
A beautiful example of this is
Cherednik's proof of
Macdonald's inner product identities
using double affine Hecke algebras.
Roughly speaking, Macdonald formulated Coxeter group analogues
of a classical combinatorial problem called the ``Dyson conjecture,''
and this paved the way for Cherednik's discovery of an algebraic
structure underlying the phenomena.  For an exposition, see~[9].

It is therefore natural to ask if
matroids can be defined in terms of the symmetric group.
If so, what do the Coxeter analogues of matroids look like?
J.~P.~S. Kung~[10, 11] seems to have been the first
to ask and answer these questions.
Later, I.~M. Gelfand and V.~V. Serganova~[8]
suggested a different definition of such analogues
(which they called ``{\it WP}-matroids'').  In both cases the work has
a strong geometric flavor, and one exciting possibility is that these
``Coxeter matroids'' may form the foundation for discrete symplectic
and orthogonal geometry, in the same way that oriented matroids form
the basis for MacPherson's theory of combinatorial differentiable
manifolds~[13].

Before such a ``combinatorial Erlanger program'' can be carried out,
however, many fundamental questions must first be answered.
A basic fact about matroids is that they admit a wide variety
of equivalent elementary axiomatizations:
independent sets, circuit elimination, basis exchange, etc.
But so far no analogous elementary axiomatizations
for {\it WP}-matroids are known.
Borovik, Gelfand and White~[2] have made some
progress in obtaining such elementary axiomatizations in the symplectic
case, but their paper does not, for example, give an independent set
axiomatization of symplectic matroids.  Indeed, for some time it was
suspected that finding such an axiomatization might be an
intractable~problem.

The main result of this paper is an independent set axiomatization for
(Gelfand-Serganova) symplectic matroids.  In addition to answering a
very natural question, this result provides one of the simplest ways
to date of explaining what a symplectic matroid is to someone with no
background in matroids or Coxeter groups.

\sec Definitions

The goal of this section is
to give the definition of a symplectic matroid.
The standard definition involves
the Bruhat order on parabolic quotients of a Coxeter group,
but in order to keep everything as simple as possible,
we will take advantage of the results in~[2]
and define symplectic matroids in a way that requires no explicit
mention of such concepts.
Readers familiar with Coxeter groups
who want the full story should see~[15].

Let $E_{\pm n}$ be the set $\{\pm1, \pm 2, \pm3, \ldots, \pm n\}$.
For brevity we will sometimes write the minus sign on top;
e.g., we will write $\bar 1$ for~$-1$.  If $w$ is a permutation
of~$E_{\pm n}$
and $B = \{ b_1, b_2, \ldots, b_k \}$ is a subset of~$E_{\pm n}$,
then we define $$wB \defeq \{wb_1, wb_2, \ldots, wb_k\}.$$
An important concept in Coxeter matroid theory is \itc{admissibility}.
If $S \subseteq E_{\pm n}$, define
$$\bar S \defeq \{-s \mid s\in S\}.$$
We say that $S$ is \itc{admissible} if $S\cap \bar S = \varemptyset$.
A permutation~$w$ of~$E_{\pm n}$ is \itc{admissible} if
$w(-x) = -wx$ for all $x\in E_{\pm n}$.
A total ordering $\prec$ of the elements of~$E_{\pm n}$ is
\itc{admissible} if there exists an admissible permutation~$w$
such that $x\prec y$ if and only if $wx < wy$.
(The reader may find it helpful
to visualize an admissible ordering
as a signed permutation~$\sigma$ of $\{1, 2, \ldots, n\}$
followed by the negative of the reversal of~$\sigma$,
e.g., $\bar2, 1, 3, \bar3, \bar1, 2$.)
If $\prec$ is an admissible total ordering of~$E_{\pm n}$,
then a map $\omega: E_{\pm n} \to \R$ is said to be a
\itc{weight function compatible with~$\prec$}
if $i \prec j$ implies $\omega(i) \le \omega(j)$.

One way of defining ordinary matroids involves the greedy algorithm
[14, section~1.8].
This is the approach we shall take to symplectic matroids.
Suppose we are given an admissible total ordering~$\prec$,
a weight function~$\omega$ compatible with $\prec$,
and a collection~$\SB$ of subsets of~$E_{\pm n}$.
Then we define the \itc{greedy solution} of~$\SB$
to be the element $B \in \SB$ that is constructed
as follows: we begin with no elements in~$B$
and then we consider each element of~$E_{\pm n}$
in turn from the largest (relative to~$\prec$) to the smallest,
adding it to~$B$ unless doing so would
make it impossible to end up with an member of~$\SB$
no matter which other elements of~$E_{\pm n}$ we subsequently add to~$B$.

For example, suppose $n=3$ and
our admissible total ordering is the standard ordering.
Let $\SB = \bigl\{ \{\bar 2, \bar 1, 3\}, \{\bar 2, 1, 3\} \bigr\}$.
We begin by putting $3$ into~$B$,
because there are certainly members of~$\SB$ containing~$3$.
We next consider~$2$, but we can't add $2$ to~$B$,
because if we do then regardless of what further numbers we add to~$B$,
we can never produce a member of~$\SB$.
(In other words, $\{2,3\}$ is not a subset of any member of~$\SB$.)
Continuing in this way,
we find that the greedy solution is $\{\bar 2, 1, 3\}$.

Finally, we say that the greedy solution~$B$ of~$\SB$ is \itc{optimal}
if $\omega(B) \ge \omega(B')$ for all $B'\in\SB$,
where as usual $\omega(B)$ denotes $\sum_{b\in B} \omega(b)$.
We can now define a symplectic matroid.

\definition
A \itc{symplectic matroid} is a pair $(E_{\pm n}, \SB)$
where $\SB$ is a nonempty family of equinumerous admissible
subsets of~$E_{\pm n}$ with the property that
for every admissible total ordering~$\prec$ of~$E_{\pm n}$
and every weight function compatible with~$\prec$,
the greedy solution of~$\SB$ is optimal.
The family~$\SB$ is called
the family of \itc{bases} of the symplectic matroid.

\remark
The equivalence of this definition of symplectic matroid with the usual
definition is the content of [2, Theorem~16].

An example of a symplectic matroid is $(E_{\pm 3}, \SB)$ where
$\SB = \{1\bar3, 2\bar3, \bar12, \bar13\}$.
Here $1\bar3$ is to be understood as shorthand for the set
$\{1, \bar 3\}$.  Note that a symplectic matroid is \itc{not}
a matroid; it is an \itc{analogue} of a matroid.
(There is a sense in which ordinary matroids may be
regarded as special cases of symplectic matroids,
but this need not concern us here.)

\sec The Main Result

If $(E_{\pm n}, \SB)$ is a symplectic matroid, we define its
\itc{family~$\I$ of independent sets} by
$$\I \defeq \{I\subseteq E_{\pm n} \mid
   I\subseteq B \;\hbox{for some $B\in \SB$}\}.$$
In the example of a symplectic matroid given in the last section,
the family of independent sets is
$\I = \{\varemptyset, 1, \bar 1, 2, 3, \bar3\} \cup \SB$.
Notice that we can recover $\SB$ from~$\I$; the members of $\SB$ are
just the maximal members of~$\I$ with respect to inclusion.
Thus, a characterization of~$\I$ could be used as
an alternative definition or axiomatization of a symplectic matroid.
This is precisely what the following theorem provides.

\theorem\mainresult
{A subset-closed family $\I$ of admissible subsets of $E_{\pm n}$ is
the family of independent sets of a symplectic matroid if and only
if it has the following property:
 
{\narrower\noindent
   If $I$ and $J$ are members of $\I$ such that $|I| < |J|$ and such
   that for all $y \in J\setminus I$, the set $\{y\} \cup I$ is not
   in~$\I$, then $I \cup J$ is inadmissible and there exists $x\notin I$
   such that both $\{x\} \cup I$ and $\{\bar x\} \cup I \setminus \bar J$
   are in~$\I$.\par}
}

We remark that part of this theorem---the part about
$I \cup J$ being inadmissible---was previously known and is essentially
Theorem~14 of~[2].
Notice incidentally that the expression
``$\{\bar x\} \cup I \setminus \bar J$'' looks ambiguous because
it is not clear whether we take the union first or subtract first,
but actually there is no ambiguity since the hypothesis
prevents $x$ from being an element of~$J$.

\proof
Unless otherwise specified, the terms ``larger'' and ``smaller''
in this proof refer to the admissible total ordering~$\prec$.
The reader should visualize such an ordering by writing out
the elements in order in a horizontal line, with the largest
element first.

\itc{Sufficiency}.
Assume that $\I$ has the stated property.
Call a maximal (with respect to inclusion) member of~$\I$
a ``basis'' of~$\I$.
All bases of~$\I$ are admissible, and the stated property of~$\I$ ensures
that all bases of~$\I$ have the same number of elements.  Let $\SB$ be the
collection of bases of~$\I$.
We now make the following claim, which we shall call~$(*)$.

\item{$(*)$}
     Let $\prec$ be an admissible ordering.  Let $I$ be a set consisting
     of the first $i$ elements of $E_{\pm n}$ that are picked up by the
     greedy algorithm (for some $i \ge 0$).
     Let $J$ be a member of~$\I$ such that $i < |J|$.  Then
     the $(i+1)$st element picked up by the greedy algorithm is
     no smaller than the smallest element of~$J$.

To see this, note first that
if there exists $y \in J\setminus I$ such that $\{y\} \cup I \in \I$,
then we are done, because then the greedy algorithm will pick
up either $y$ or something larger than~$y$, and $y$ is trivially no
smaller than the smallest element of~$J$.
Otherwise, since $\I$ has the stated property, there exists $x\notin I$
such that $\{x\} \cup I$ and $\{\bar x\} \cup I \setminus \bar J$
are both in~$\I$.  Moreover, $I \cup J$ is
inadmissible, but each of $I$ and~$J$ is admissible,
so there exists $z\in I$ such that $\bar z\in J$.
Choose the largest such~$z$.
Then by the maximality of~$z$, the set $S$ of elements of $I$ that
are larger than~$z$ is a subset of $I \setminus \bar J$,
and therefore both $S \cup \{x\}$ and $S \cup \{\bar x\}$ are in~$\I$.
Now $\bar x \notin I$ (since $\{x\} \cup I$ is in $\I$ and is therefore
admissible) and also $x \notin I$, so neither $x$ nor $\bar x$ can be
larger than~$z$---otherwise, since both $x$ and~$\bar x$
are ``compatible'' with~$S$, the greedy algorithm
would have picked one of them
(or some other element not in~$I$ that is even larger), and it didn't.
It follows that $z$ appears before the ``halfway point''
(the point between the $n$th and the $(n+1)$st elements in the ordering),
and that $x$ and $\bar x$ both appear after~$z$ but before~$\bar z$.
Then the $(i+1)$st element picked up by the greedy algorithm
must be no smaller than~$x$, which is no smaller then~$\bar z$,
which in turn is no smaller than the smallest element of~$J$,
since $\bar z \in J$.  This proves~$(*)$.

Now let $\prec$ be an admissible ordering
and let $\omega$ be some weight function
compatible with~$\prec$.
Let $B$ be the basis of~$\I$ chosen by the greedy algorithm.
We want to show that $B$ is optimal, so let $B'$ be another basis.
We claim that for all $i > 0$,
the $i$th element of~$B$ is no smaller than
the $i$th element of~$B'$.  For, given~$i$, let $I$ be the set consisting
of the largest $i-1$ elements of~$A$
and let $J$ be the set consisting of the largest
$i$ elements of~$B'$.
Then $|I| < |J|$, so by~$(*)$, the $i$th element picked up by
the greedy algorithm (i.e., the $i$th element of~$B$)
is no smaller than the smallest element
(i.e., the $i$th element) of~$J$.
This proves the claim,
which in turn shows that for all~$i$, the weight of the $i$th element
of~$B$ is greater than or equal to
the weight of the $i$th element of~$B'$, so indeed $B$ is optimal.

\itc{Necessity}.
Suppose that $\I$ is the family of independent sets of
a symplectic matroid, and let $I$ and~$J$ be members of~$\I$
such that $|I| < |J|$
and such that for all $y \in J\setminus I$,
$\{y\} \cup I \notin \I$.
Then, as already mentioned, [2, Theorem~14]
implies that $I \cup J$ is inadmissible.
The set $I \cup J$ may be partitioned into
four disjoint sets $W$, $Y$, $Z$, and $\bar Z$,
where $W$, $Y$, and~$Z$ are defined as follows:
$$\eqalign{ W &= I \setminus \bar J\cr
            Y &= J \setminus (I \cup \bar I)\cr
            Z &= I \cap \bar J\cr}$$
In words, $Z$ is the subset of~$I$ whose negatives are in~$J$,
$W$ is the rest of~$I$,
and $Y$ is what remains in~$J$ after $W$, $Z$, and~$\bar Z$ are removed.

Now let
$X = E_{\pm n} \setminus (W \cup \bar W \cup
          Y \cup \bar Y \cup Z \cup \bar Z)$.
Define a ``half'' of~$X$ to be
a maximal (with respect to inclusion) admissible subset of~$X$.
Clearly, if $H$ is a half of~$X$,
then $H$ and $\bar H$ partition~$X$ into two disjoint sets
and $|H| = |\bar H|$.
Define a ``\itc{WXYZ} ordering'' to be an admissible ordering
in which the elements of~$W$ come first,
then the elements of some half~$H$ of~$X$,
then the elements of~$Y$, and then the elements of~$Z$.
(This gives us half of~$E_{\pm n}$,
so the ordering of
the rest of~$E_{\pm n}$---namely
$\bar Z \bar Y \bar H \bar W$---is determined.)

Now suppose we are given a \itc{WXYZ} ordering with the weight function that
equals one on $W$, $H$, $Y$, $Z$, and~$\bar Z$ and equals zero after that.
The greedy algorithm will begin by picking up the elements of~$W$.
We claim that some element of $H \cup Y$ must be picked up after that.
For if not, the algorithm will pick up~$Z$,
since these are just the remaining elements of~$I$.
Then it will skip over~$\bar Z$.
This implies that the weight of the basis chosen will be~$|I|$,
but $J$ is contained in $W \cup Y \cup \bar Z$
so the weight of~$J$ is $|J| > |I|$, a contradiction.

The argument just given applies
regardless of how the half~$H$ of~$X$ is chosen.
Therefore the following set~$S$ is nonempty:
$$ S = \bigl\{x \in X \mid \hbox{$\{x\} \cup W \in \I$
                             and $\{\bar x\} \cup W \in \I$} \bigr\}
 \cup  \bigl\{y \in Y \mid \{y\} \cup W \in F\bigr\}.$$
(For if not, we could choose a half $H$ of~$X$ such that for
\itc{all} $x\in H$,
$\{x\} \cup W$ would not be in~$\I$,
and this would cause trouble for the greedy
algorithm as just explained.)  Now construct an admissible ordering~$\prec$
as follows.  Begin with a \itc{WXYZ} ordering that minimizes the number of
$x \in H$ such that $\{x\} \cup W \in \I$.
Then reposition every element in
$(H \cup Y) \cap S$ so that they now come after~$Z$ (but before~$\bar Z$).
Finally, reposition the ``mirror images'' of the elements just moved
to restore admissibility.  For example, if the \itc{WXYZ} ordering were
\def\ds{\displaystyle}
$$\underbrace{a\quad b\strut}_{\ds W} \quad
  \underbrace{c\quad d\strut}_{\ds H} \quad
  \underbrace{e\quad f\strut}_{\ds Y} \;
  \underbrace{g\strut}_{\ds Z} \;
  \bar g \quad \bar f \quad \bar e \quad \bar d \quad
  \bar c \quad \bar b \quad \bar a$$
and $d$ and~$e$ were in~$S$ but $c$ and~$f$ were not,
then $\prec$ would be given by
$$ a \succ b \succ c \succ f \succ g \succ d \succ e \succ
  \bar e \succ \bar d \succ \bar g \succ \bar f \succ \bar c
  \succ \bar b \succ \bar a.$$

Observe that by the minimality in the choice of~$H$,
the elements $x \in H \cup Y$ that are not repositioned
have the property that $\{x\} \cup W \notin H$.
Now give every element up to the end of~$\bar Z$ weight one and give
the rest of the elements weight zero.  The greedy algorithm applied
to this ordering will pick up the elements of~$W$, and will skip over
the elements of~$H \cup Y$.  Then it will pick up the elements of~$Z$,
since (as before) these are just the remaining elements of~$I$.
Now, as before, $J$ has greater weight than~$I$,
so the greedy algorithm must pick up another
element before it reaches the end of~$\bar Z$.
It cannot pick up any element of~$\bar Z$,
so it must pick up one of the repositioned elements
($d$, $e$, $\bar e$, or~$\bar d$ in the example above).
Let $x$ be the first element so picked up.
If $x \in X$, then we see that it satisfies the desired conditions
(that both $\{x\} \cup I$ and $\{\bar x\} \cup I \setminus \bar J$
are in~$\I$).  Otherwise, $x$ cannot be in~$Y$,
because $Y \subseteq J$ and for no $x \in J\setminus I$
can we have $\{x\} \cup I \in \I$.
So $x \in \bar Y$.
In particular, $x \in \bar J$,
so $\{\bar x\} \cup I \setminus \bar J = I$,
which is trivially in~$\I$.  This completes the proof.\qed

\sec What Next?

It would be nice to find an independent set axiomatization
for \itc{all} Coxeter matroids, not just symplectic ones.
We might also hope to use \mainresult\ to obtain circuit
elimination and basis exchange axioms for Coxeter matroids,
since these axioms are closely related to independent set
axioms in the ordinary matroid case.

It should also be fruitful to determine the precise connections
between Gelfand-Serganova Coxeter matroids and
all the other generalizations of matroids that exist in the literature.
Following are some observations
that were revealed by a quick literature search.

There is one special case of a symplectic matroid
that has been rediscovered independently several times.
It goes by different names: ``Lagrangian matroid,'' ``symmetric matroid''~[3],
``$\Delta$-matroid''~[3], and ``pseudomatroid''~[6].
All these concepts are equivalent,
and Gelfand-Serganova symplectic matroids
are strictly more general than all of them, as noted in~[2].
In addition, there exists something called a ``metroid''~[7]
which is almost equivalent to a $\Delta$-matroid,
but technically it is a special case:
metroids are $\Delta$-matroids that include the empty set as a feasible set.
This is proved in~[4].
Incidentally, the Mathematical Review 89a:05046 of~[3]
remarks that ``symmetric matroid'' is also used in the literature
to refer to something completely different,
but I have not been able to track down any instances of this other usage.

A concept that is earlier than any of the above
is that of a ``bimatroid'' [10,~11].
In~[7] it is shown that a bimatroid is a special case of a metroid.
In~[11], two concepts that are related to bimatroids are discussed:
``orthogonal matroids'' and ``Pfaffian structures.''
Orthogonal matroids are special cases of bimatroids
and hence (confusingly) are special cases of
Gelfand-Serganova \itc{symplectic} matroids.
Gelfand-Serganova orthogonal matroids (i.e., the case $W=D_n$)
may also be viewed as special cases of
Gelfand-Serganova symplectic matroids,
but it is not immediately clear whether there is any direct connection
between orthogonal matroids in the sense of~[11]
and orthogonal matroids in the sense of Gelfand-Serganova.
To add to the confusion,
sometimes Pfaffian structures are referred to as ``symplectic matroids''
because they are indeed symplectic analogues of matroids [11,~12]
but it is not immediately clear what the precise relationship
between them and the other concepts mentioned above is.
One can get a Pfaffian structure out of a bimatroid,
but they do not seem to be strictly equivalent,
and thus a Pfaffian structure does not seem to be
a special case of (say) a $\Delta$-matroid.

Finally, two other concepts that might be related to Coxeter matroids
are ``Coxeteroids''~[1] and ``multimatroids''~[5].
The definition of a Coxeteroid is motivated by the observation that
matroids and Coxeter groups both satisfy an ``exchange condition.''
A multimatroid is a certain generalization of a $\Delta$-matroid.
In neither case is the exact connection with Coxeter matroids obvious.

\sec Acknowledgments

This work was supported in part by a National Science Foundation
Postdoctoral Fellowship.  Part of the work was done while I was a
general member of the Mathematical Sciences Research Institute.
I thank Neil White for introducing me to this problem and for
helpful discussions and encouragement.

\sec References

\parskip=5pt plus1pt minus1pt

{\smrm
\advance\baselineskip by -1.5pt
\def\smitc#1{{\smit #1\/}}

\item{1.} {\smsc A. Bj\"orner},
On matroids, groups, and exchange languages, \smitc{in}
``Matroid Theory,'' ed.\ L.~Lov\'asz and A.~Recski,
North Holland, New York, 1985.

\item{2.} {\smsc A. V. Borovik, I. M. Gelfand, and N. L. White},
Symplectic matroids, \smitc{J.~Alg.\ Combin.,} accepted for publication.

\item{3.} {\smsc A. Bouchet}, Greedy algorithm and symmetric matroids,
\smitc{Math.\ Prog.}\ {\smbf 38} (1987), 147--159.

\item{4.} {\smsc A. Bouchet, A. W. M. Dress and T. F. Havel},
$\Delta$-matroids and metroids, \smitc{Advances in Math.}\
{\smbf 91} (1992), 136--142.

\item{5.} {\smsc A. Bouchet}, Multimatroids~I.~Coverings by independent
sets, \smitc{SIAM J.~Discrete Math.}\ {\smbf 10} (1997), 626--646.

\item{6.} {\smsc R. Chandrasekaran and S. N. Kabadi}, Pseudomatroids,
\smitc{Discrete Math.}\ {\smbf 71} (1988), 205--217.

\item{7.} {\smsc A. Dress and T. F. Havel}, Some combinatorial properties
of discriminants in metric vector spaces, \smitc{Advances in Math.}\
{\smbf 62} (1986), 285--312.

\item{8.} {\smsc I. M. Gelfand and V. V. Serganova},
Combinatorial geometries and
torus strata on homogeneous compact manifolds, \smitc{Russ.\ Math.\ Surv.}\
{\smbf 42} (1987), 133--168.

\item{9.} {\smsc A. A. Kirillov, Jr.,} Lectures on affine Hecke algebras and
Macdonald's conjectures, \smitc{Bull.\ Amer.\ Math.\ Soc.}\ {\smbf 34}
(1997), 251--292.

\item{10.} {\smsc J. P. S. Kung},
An Erlanger program for combinatorial geometries,
Ph.D. thesis (1978), Massachusetts Institute of Technology, Cambridge,~MA.

\item{11.} {\smsc J. P. S. Kung,} Bimatroids and invariants,
\smitc{Advances in Math.}\ {\smbf 30} (1978), 238--249.

\item{12.} {\smsc J. P. S. Kung}, Pfaffian structures and critical
problems in finite symplectic spaces, \smitc{Ann.\ Combin.}\
{\smbf 1} (1997), 159--172.

\item{13.} {\smsc R. D. MacPherson}, Combinatorial differential manifolds,
in ``Topological Methods in Modern Mathematics: A Symposium in Honor of
John Milnor's Sixtieth Birthday,'' ed.~L.~R. Goldberg and A.~V. Phillips,
Publish or Perish, Inc., Houston, 1993, 203--221.

\item{14.} {\smsc J. G. Oxley}, ``Matroid Theory,'' Oxford University Press,
New York, 1992.

\item{15.} {\smsc N. L. White}, The Coxeter matroids of Gelfand et al.,
\smitc{Contemp.\ Math.}\ {\smbf 197} (1996), 401--409.

}
\bye